\documentclass{article}
\usepackage{amssymb}
\usepackage{amsmath}

\input{tcilatex}
\begin{document}

\begin{center}
{\LARGE Nonlinear differential algorithm to compute all the zeros of a
generic polynomial}

\bigskip

\textbf{Francesco Calogero}

Physics Department, University of Rome "La Sapienza", 00185 Rome, Italy

Istituto Nazionale di Fisica Nucleare, Sezione di Roma

\bigskip

\textit{Abstract}
\end{center}

A simple algorithm to compute all the zeros of a generic polynomial is
proposed.

\bigskip

\section{Introduction}

\textbf{Notation}. Hereafter, for definiteness, we always refer to \textit{%
monic} polynomials of arbitrary order $N$ ($N\geq 2$),%
\begin{equation}
P_{N}\left( z;\vec{c},\underline{x}\right) =z^{N}+\sum_{m=1}^{N}\left(
c_{m}~z^{N-m}\right) =\tprod\limits_{n=1}^{N}\left( z-x_{n}\right) ~;
\label{Pol}
\end{equation}%
the \textit{complex} variable $z$ is the argument of the polynomial, indices
such as $n$, $m$, $\ell $ run from $1$ to $N$ (unless otherwise indicated,
see below), the $N$-vector $\vec{c}$ has the $N$ \textit{coefficients }$%
c_{m} $ of the polynomial (\ref{Pol}) as its $N$ components, $\underline{x}$
is the \textit{unordered} set of the $N$ \textit{zeros} $x_{n}$ of the
polynomial (\ref{Pol}), and we assume all variables to be \textit{complex}
(unless otherwise explicitly indicated, see below)\textit{.} We call \textit{%
generic} any polynomial the \textit{coefficients} and \textit{zeros} of
which are \textit{generic complex} numbers, and in particular feature 
\textit{zeros} which are \textit{all different among themselves, }$x_{n}\neq
x_{m}$ if $n\neq m$. Note that the notation $P_{N}\left( z;\vec{c},%
\underline{x}\right) $ is somewhat redundant, since this monic polynomial
can be identified by assigning \textit{either} its $N$ coefficients \textit{%
or} its $N$ zeros; indeed the $N$ coefficients $c_{m}$ can be expressed in
terms of the $N$ zeros $x_{n}$ via the following standard formula 
\begin{subequations}
\label{cm}
\begin{equation}
c_{m}=\left( -1\right) ^{m}~\sum_{n_{1}>n_{2}>...>n_{m}=1}^{N}\left(
x_{n_{1}}~x_{n_{2}}\cdot \cdot \cdot x_{n_{m}}\right) ~,
\end{equation}%
so that%
\begin{equation}
c_{1}=-\left( x_{1}+x_{2}+...+x_{N}\right) ~,
\end{equation}%
\begin{eqnarray}
&&c_{2}=\left( x_{1}~x_{2}+x_{1}~x_{3}+...+x_{1}~x_{N}\right)  \notag \\
&&+\left( x_{2}~x_{3}+x_{2}~x_{4}+...+x_{2}~x_{N}\right) +...  \notag \\
&&+\left( x_{N-2}~x_{N-1}+x_{N-2}~x_{N}\right) +x_{N-1}~x_{N}~,
\end{eqnarray}%
and so on. On the other hand, while the assignment of the $N$ \textit{%
coefficients} $c_{m}$ determines the $N$ \textit{zeros }$x_{n}$---uniquely,
up to permutations---of course \textit{explicit} formulas to accomplish
generally this task \textit{only} exist for $N\leq 4$. $\blacksquare $

The investigation of the properties---and of techniques for the numerical
computation---of the $N$ \textit{zeros} $x_{n}$ of a polynomial of degree $N$
defined via the assignment of its $N$ \textit{coefficients} $c_{m}$ (see (%
\ref{Pol})) is a problem that has engaged mathematicians since time
immemorial. In this paper a simple nonlinear differential algorithm suitable
to compute numerically \textit{all} the $N$ zeros of a \textit{generic}
polynomial of arbitrary degree $N$ is described; I was unable to find a
previous description of this algorithm in the literature, but I am aware
that my search has not been---indeed, it could not have been---quite
complete. This algorithm is described in the following Section 2 and proven
in Section 3.

\bigskip

\section{Results}

It is now convenient to introduce an additional independent variable $t$,
which is hereafter assumed to be \textit{real} and might be interpreted as 
\textit{time}. Hence the above notation is now extended by writing, in
addition to (\ref{Pol}), the analogous formula 
\end{subequations}
\begin{equation}
p_{N}\left( z;\vec{\gamma}\left( t\right) ,\underline{y}\left( t\right)
\right) =z^{N}+\sum_{m=1}^{N}\left[ \gamma _{m}\left( t\right) ~z^{N-m}%
\right] =\tprod\limits_{n=1}^{N}\left[ z-y_{n}\left( t\right) \right] ~,
\label{pol}
\end{equation}%
to which notational comments quite analogous to those reported above apply.

There holds then the following

\textbf{Proposition}. Consider the following system of $N$ nonlinear
first-order differential equations satisfied by the $N$ \textit{zeros} $%
y_{n}\left( t\right) $ of the polynomial (\ref{pol}): 
\begin{subequations}
\label{ydotgamma}
\begin{equation}
\dot{y}_{n}\left( t\right) =-\left\{ \tprod\limits_{\ell =1,~\ell \neq n}^{N}%
\left[ y_{n}\left( t\right) -y_{\ell }\left( t\right) \right] ^{-1}\right\}
~\sum_{m=1}^{N}\left\{ \left[ c_{m}-\gamma _{m}\left( 0\right) \right] ~%
\left[ y_{n}\left( t\right) \right] ^{N-m}\right\} ~.  \label{ydot}
\end{equation}%
Here and below a superimposed dot denote a $t$-differentiation, while the 
\textit{coefficients} $c_{m}$ are those of the polynomial $P_{N}\left( z;%
\vec{c},\underline{x}\right) $, see (\ref{Pol}), the \textit{zeros} $x_{n}$
of which we seek, and $\gamma _{m}\left( 0\right) $ are the $N$ \textit{%
coefficients} of the polynomial $p_{N}\left( z;\vec{\gamma}\left( t\right) ,%
\underline{y}\left( t\right) \right) $, see (\ref{pol}), at $t=0,$ hence
they are related to the "initial" values $y_{n}\left( 0\right) $ of the 
\textit{zeros} of this polynomial by the formula (analogous to (\ref{cm}))%
\begin{equation}
\gamma _{m}\left( 0\right) =\left( -1\right)
^{m}~\sum_{n_{1}>n_{2}>...>n_{m}=1}^{N}\left[ y_{n_{1}}\left( 0\right)
~y_{n_{2}}\left( 0\right) \cdot \cdot \cdot y_{n_{m}}\left( 0\right) \right]
~.  \label{gammamzero}
\end{equation}

Then 
\begin{equation}
x_{n}=y_{n}\left( 1\right) ~.~~~\blacksquare  \label{xyn}
\end{equation}

It is thus seen that the \textit{zeros} $x_{n}$ of the polynomial $%
P_{N}\left( z;\vec{c},\underline{x}\right) ,$ see (\ref{Pol}), can be
computed---once the $N$ coefficients $c_{m}$ of this polynomial have been
assigned---via the following procedure. \textit{Step one}: choose (\textit{%
arbitrarily}!) $N$ \textit{complex} numbers $y_{n}\left( 0\right) $. \textit{%
Step two}: compute, via the formulas (\ref{gammamzero}), the $N$ quantities $%
\gamma _{m}\left( 0\right) $. \textit{Step three}: integrate (numerically)
the system of differential equations (\ref{ydot}) from $t=0$ to $t=1$,
starting from the $N$ initial data $y_{n}\left( 0\right) $, getting thereby
the $N$ values $y_{n}\left( 1\right) $, which give the sought result, see (%
\ref{xyn}).

Will this procedure always work?\ The only possible snag is that the
solution $\vec{y}\left( t\right) $ of the "dynamical system" (\ref{ydot})
run into a singularity during its evolution from $t=0$ to $t=1.$ The only
mechanism whereby this might occur is because during this evolution two
different coordinates $y_{n}\left( t\right) $ might coincide, $y_{\ell
}\left( t\right) =y_{n}\left( t\right) $ for $\ell \neq n$, at some value of
the \textit{real} variable $t$ in the interval $0<t<1$, causing the
right-hand side of (\ref{ydot}) to blow up. This "collision" might indeed
happen, but it is \textit{not} a \textit{generic} phenomenon: hence it will
be enough to change the assignment of the (arbitrary!) initial data $%
y_{n}\left( 0\right) $ to avoid this difficulty; note however that this
suggests that to apply this method it will be advisable to always start with 
\textit{complex} initial data $y_{n}\left( 0\right) $, even in the case of 
\textit{real} polynomials with \textit{real} zeros. And note that by
performing the numerical integration of the differential equations (\ref%
{ydot}) with different initial data $\vec{y}\left( 0\right) $ provides
moreover a possibility to assess the \textit{numerical accuracy} of the
computation, by comparing the results obtained starting from different sets
of initial data.

\textbf{Remark}. It is plain that this procedure will work more efficiently
the closer the, arbitrarily chosen, initial values $y_{n}\left( 0\right) $
are to the $N$ \textit{zeros} $x_{n}$ the values of which one is trying to
compute; indeed if the $N$ initial values $y_{n}\left( 0\right) $ happened
to \textit{coincide} with the $N$ zeros $x_{n}$, $y_{n}\left( 0\right)
=x_{n} $, this would imply $\gamma _{m}\left( 0\right) =c_{m}$ (compare (\ref%
{cm}) with (\ref{gammamzero})) hence the right-hand side of the differential
equations (\ref{ydot}) would vanish identically, entailing $\dot{y}_{n}=0$
hence $y_{n}\left( 1\right) =y_{n}\left( 0\right) =x_{n}$, consistently with
(\ref{xyn}).

Let us also emphasize that the dependence (via (\ref{gammamzero})) of the
right-hand sides of the differential equations (\ref{ydot}) upon the initial
values $y_{n}\left( 0\right) $ of the dependent variables $y_{n}\left(
t\right) $ implies that these differential equations are rather Differential
Functional Equations than Ordinary Differential Equations; but this fact has
hardly any relevance on \textit{step three} of the procedure, see above. $%
\blacksquare $

A comparison of the actual effectiveness of this technique with that of
other methods to compute \textit{all} the $N$ zeros of a \textit{generic}
polynomial of arbitrary degree $N$ is beyond the scope of this short
communication, and in any case it is a task to be rather pursued by
specialists in numerical analysis if they consider it worthy of their
attention.

\bigskip

\section{Proof}

The proof of the above \textbf{Proposition} is actually quite easy (raising
thereby some doubts on the novelty of this finding). The starting point is
the \textit{identity} 
\end{subequations}
\begin{equation}
\dot{y}_{n}\left( t\right) =-\left\{ \tprod\limits_{\ell =1,~\ell \neq n}^{N}%
\left[ y_{n}\left( t\right) -y_{\ell }\left( t\right) \right] ^{-1}\right\}
~\sum_{m=1}^{N}\left\{ \dot{\gamma}_{m}\left( t\right) ~\left[ y_{n}\left(
t\right) \right] ^{N-m}\right\} ~,  \label{Iden}
\end{equation}%
valid for any $t$-dependent polynomial with \textit{zeros} $y_{n}\left(
t\right) $ and \textit{coefficients} $\gamma _{m}\left( t\right) ,$ see (\ref%
{pol}); for a proof of this formula see \cite{C2015}. Now make the
assignment 
\begin{subequations}
\begin{equation}
\gamma _{m}\left( t\right) =\gamma _{m}\left( 0\right) +\left[ c_{m}-\gamma
_{m}\left( 0\right) \right] ~t~,  \label{gammat1}
\end{equation}%
consistent with the initial (arbitrary) assignment at $t=0$ and clearly
implying%
\begin{equation}
\dot{\gamma}_{m}\left( t\right) =c_{m}-\gamma _{m}\left( 0\right) ~,
\label{gammadot}
\end{equation}%
\begin{equation}
\gamma _{m}\left( 1\right) =c_{m}~.  \label{gammam1}
\end{equation}%
The insertion of the first of these two formulas, (\ref{gammadot}), in (\ref%
{Iden}) yields (\ref{ydot}); while the second, (\ref{gammam1}), implies
that, at $t=1$, the polynomial $p_{N}\left( z;\vec{\gamma}\left( t\right) ,%
\underline{y}\left( t\right) \right) $, see (\ref{pol}), coincides with the
polynomial $P_{N}\left( z;\vec{c},\underline{x}\right) ,$ see (\ref{Pol}),
hence the validity of (\ref{xyn}). Q. E. D.

\end{subequations}

\end{document}